\documentclass[12pt,intlim]{amsart}

\usepackage{amssymb,amsmath}

\usepackage[arrow,matrix,curve]{xy}

\newdir{ >}{{}*!/-5pt/\dir{>}}

\newcommand{\nc}{\newcommand}
\nc{\al}{\alpha}
\nc{\bt}{\beta}
\nc{\gm}{\gamma}
\nc{\dl}{\delta}
\nc{\sg}{\sigma}
\nc{\Sg}{\Sigma}
\nc{\vf}{\varphi}
\nc{\ve}{\varepsilon}
\nc{\ol}{\overline}
\nc{\lb}{\lambda}
\nc{\Lb}{\Lambda}
\nc{\vkp}{\varkappa}
\nc{\ul}{\underline}
\nc{\pa}{\partial}
\nc{\bZ}{\Bbb{Z}}
\nc{\tm}{\times}
\nc{\sbs}{\subset}
\nc{\lra}{\longrightarrow}
\nc{\all}{\allowdisplaybreaks}
\nc{\Ker}{\operatorname{Ker}}
\nc{\cl}{\operatorname{cl}}

\begin{document}

\title[ON EXACTNESS OF LONG SEQUENCES]
{ON EXACTNESS OF LONG SEQUENCES OF HOMOLOGY SEMIMODULES}

\author{Alex Patchkoria}

\begin{abstract}
We investigate exactness of long sequences of homology semimodules
associated to Schreier short exact sequences of chain complexes of
semimodules.
\end{abstract}

\thanks{Partially supported by INTAS grant 00-566}



\maketitle

In [4], to define homology and cohomology monoids of presimplicial
semimodules, we
 introduced a chain complex of semimodules (in particular, abelian monoids),
its homology semimodules, and a $\pm$-morphism between chain
complexes of semimodules. Next, in [5], we introduced a morphism
between chain complexes of cancellative semimodules,
 defined a chain homotopy of morphisms and studied its basic properties.
 In this paper we investigate exactness of long sequences of homology semimodules
 associated to Schreier short exact sequences of chain complexes of semimodules.

The paper is divided into two sections. Section 1 contains the
preliminaries. The main  results are presented in Section 2.

\section{Preliminaries}

Recall [1] that a semiring $\Lb=(\Lb,+\,,0,\cdot\,,1)$ is an
algebraic  structure in which  $(\Lb,+\,,0)$   is an abelian
monoid, $(\Lb,\cdot\,,1)$ a monoid, and
\begin{align*}
\lb\cdot(\lb'+\lb'')&=\lb\cdot\lb'+\lb\cdot\lb'',\\
(\lb'+\lb'')\cdot\lb&=\lb'\cdot\lb+\lb''\cdot\lb,\\
\lb\cdot 0=0\cdot\lb&=0,
\end{align*}
for all $\lb,\lb',\lb''\in\Lb$. An abelian monoid $A=(A,+\,,0)$ together with a map
$\Lb\tm A\lra A$, written as $(\lb,a)\mapsto\lb a$, is called a (left) $\Lb$-semimodule if
\begin{align*}
\lb(a+a')&=\lb a+\lb a',\\
(\lb+\lb')a&=\lb a+\lb' a,\\
(\lb\cdot \lb')a&=\lb(\lb' a),\\
1a=a,&\quad 0a=0,
\end{align*}
for all $\lb,\lb'\in \Lb$ and $a,a'\in A$. It immediately follows that $\lb 0=0$ for
any $\lb\in\Lb$.

A map $f:A\lra B$ between $\Lb$-semimodules $A$ and $B$ is called
a $\Lb$-homomorphism if $f(a+a')=f(a)+f(a')$ and $f(\lb a)=\lb
f(a)$, for all $a,a'\in A$ and  $\lb\in\Lb$. It is obvious that
any $\Lb$-homomorphism carries 0 into 0. A $\Lb$-subsemimodule $A$
of a   $\Lb$-semimodule $B$ is a subsemigroup of $(B,+)$ such that
$\lb a\in A$  for all    $a\in A$ and $\lb\in\Lb$. Clearly $0\in
A$. The quotient $\Lb$-semimodule $B/A$ is    defined as the
quotient $\Lb$-semimodule of $B$ by the smallest congruence on the
$\Lb$-semimodule $B$ some class of which contains $A$. Denote the
congruence class of $b\in B$ by $[b]$. Then  $[b_1]=[b_2]$ if and
only if $a_1+b_1=a_2+b_2$ for some $a_1,a_2\in A$.

Let $N$ be the semiring of nonnegative integers. An $N$-semimodule $A$ is simply an abelian
monoid, and an $N$-homomorphism $f:A\lra B$ is just a homomorphism of abelian monoids, and
 $A$ is an $N$-subsemimodule of an $N$-semimodule $B$ if and only if $A$ is a submonoid of
  the monoid $(B,+\,,0)$.

Next recall that the group completion of an abelian monoid $M$ can
be constructed in the following way.
 Define an equivalence relation $\sim$ on $M\tm M$ as follows:
$$(u,v)\sim(x,y)\Leftrightarrow u+y+z=v+x+z\quad\text{for some}\quad z\in M.$$

Let $[u,v]$ denote the equivalence class of $(u,v)$. The quotient set
$(M\tm M)/\!\!\sim$ with the addition $[x_1,y_1]+[x_2,y_2]=[x_1+x_2,y_1+y_2]$ is an abelian
 group
$(0=[x,x],\;-[x,y]=[y,x])$. This group, denoted by $K(M)$, is the group completion of $M$,
 and $k_M:M\lra K(M)$ defined by $k_M(x)=[x,0]$ is the  canonical homomorphism. If $M$ is
 a semiring, then the multiplication $[x_1,y_1]\cdot[x_2,y_2]=[x_1x_2+y_1y_2,x_1y_2+y_1x_2]$
 converts  $K(M)$ into the ring completion of the semiring $M$, and $k_M$ into the canonical
  semiring homomorphism. Now assume that $A$ is a $\Lb$-semimodule. Then $K(A,+,0)$ with the
   multiplication $[\lb_1,\lb_2][a_1,a_2]=[\lb_1a_1+\lb_2a_2,\lb_1a_2+\lb_2a_1]$,
$\lb_1,\lb_2\in\Lb$, $a_1,a_2\in A$,  becomes a $K(\Lb)$-module. This $K(\Lb)$-module,
 denoted by $K(A)$, is the $K(\Lb)$-module completion of the $\Lb$-semimodule $A$, and
$k_A=k_{(A,+,0)}$ is the canonical $\Lb$-homomorphism. Clearly,
$K(A)$ is in fact an additive functor: for any homomorphism
$f:A\lra B$ of $\Lb$-semimodules, $K(f):K(A)\lra K(B)$ defined by
$K(f)([a_1,a_2])=[f(a_1),f(a_2)]$ is a $K(\Lb)$-homomorphism.

A $\Lb$-semimodule $A$ is said to be cancellative if whenever  $a+a'=a+a''$, $a,a'a''\in A$,
 one has $a'=a''$. Obviously, $A$ is cancellative if and only if the canonical
 $\Lb$-homomorphism $k_A:A\lra K(A)$ is injective. Also note that $A$ is a cancellative
 $\Lb$-semimodule if and only if $A$ is a cancellative $C(\Lb)$-semimodule, where $C(\Lb)$
  denotes the largest additively cancellative homomorphic image of the semiring $\Lb$. ($C(\Lb)=\Lb/\!\!\sim$, \  $\lb_1\sim\lb_2$,  $\lb_1,\lb_2\in\Lb\Leftrightarrow\lb+\lb_1=\lb+\lb_2$, $\lb\in\Lb$. \ $\cl_\sim(\lb_1)+\cl_\sim(\lb_2)=\cl_\sim(\lb_1+\lb_2)$,\  \
$\cl_\sim(\lb_1)\cdot\cl_\sim(\lb_2)=\cl_\sim(\lb_1\cdot\lb_2)$.)

A $\Lb$-semimodule $A$ is called a $\Lb$-module if $(A,+,0)$ is an
abelian group. One can easily see that  $A$ is a $\Lb$-module  if
and only if $A$ is a $K(\Lb)$-module. Consequently, if $A$ is a
$\Lb$-module, then $K(A)=A$ and $k_A=1_A$.

\

\

{\bf 1.1.}  {\bf Definition ([7, 2, 8, 3]).} A sequence
$\xymatrix{E:A\ar@{ >-{>}}[r]^-\lb&B\ar@{ ->>}[r]^-\tau&C}$ of $\Lb$-semimodules and
$\Lb$-homomorphisms is called a Schreier extension of $A$ by $C$ (some authors would
say ``$C$ by $A$'') if the following conditions hold:
\begin{enumerate}
\item \ $\lb$ is injective, $\tau$ is surjective, and $\lb(A)=\Ker(\tau)$.
\item \  For any $c\in C$, $\tau^{-1}(c)$ contains an element $u_c$ such that for any
$b\in\tau^{-1}(c)$ there exists a unique element $a\in A$ with $b=\lb(a)+u_c$.
\end{enumerate}
The elements $u_c$, $c\in C$, are called representatives of the extension
 $\xymatrix{E:A\ar@{ >-{>}}[r]^-\lb&B\ar@{ ->>}[r]^-\tau&C}$.

 \

 \

The following four properties of Schreier  extensions of $\Lb$-semimodules are easy to verify.
\vskip+2mm

{\bf  1.2.} \ Let $\xymatrix{E:G\ar@{ >-{>}}[r]&B\ar@{ ->>}[r]&C}$
be a Schreier extension with $G$ a $\Lb$-module. Then any $b\in B$
is a representative of the extension $E$. \vskip+2mm

{\bf  1.3.} \ Let $\xymatrix{E:A\ar@{ >-{>}}[r]^-\lb&B\ar@{ ->>}[r]^-\tau&C}$ be a Schreier
extension with $A$ a cancellative $\Lb$-semimodule.  If
$\lb(a)+b_1=\lb(a)+b_2$, $a\in A$, $b_1,b_2\in B$, then $b_1=b_2$.
\vskip+2mm

{\bf  1.4.} \ If $\xymatrix{E:A\ar@{ >-{>}}[r]&B\ar@{ ->>}[r]&C}$ is a Schreier extension of
 $\Lb$-semi\-mo\-du\-les, then  $B$ is cancellative if and only if $A$ and $C$ are both
 cancellative.
\vskip+2mm

{\bf  1.5.} \ If $\xymatrix{E:A\ar@{ >-{>}}[r]^-\lb&B\ar@{ ->>}[r]^-\tau&C}$ is a Schreier
 extension of $\Lb$-semi\-mo\-du\-les, then
$\xymatrix{K(E):0\to K(A)\ar[r]^-{K(\lb)}&K(B)\ar[r]^-{K(\tau)}&K(C)\to 0}$ is a short exact
sequence of  $K(\Lb)$-modules.
\vskip+4mm

{\bf  1.6.} \  A  homomorphism $\vf:A\lra B$ of  $\Lb$-semimodules is said to be {\em normal}
 (or kernel-regular in the sense of [9]) if whenever $\vf(a_1)=\vf(a_2)$, $a_1,a_2\in A$, one
  has
$\kappa_1+a_1=\kappa_2+a_2$ for some $\kappa_1,\kappa_2\in\Ker(\vf)$. It is easy to see that
$\vf$ is normal if and only if $\vf:A\lra\vf(A)$ is a cokernel of the inclusion $\Ker(\vf)
\hookrightarrow A$ (i.e., $\vf:A\lra\vf(A)$ is a normal $\Lb$-epimorphism).
\vskip+2mm

{\bf  1.7.} \ Any $\Lb$-homomorphism $\vf:G\lra B$ with $G$ a $\Lb$-module is evidently
normal. Moreover, any $\Lb$-homomorphism $\vf:A\lra B$ with $\vf(A)$ a $\Lb$-module is normal. Consequently, if a sequence of $\Lb$-se\-mimodules and $\Lb$-homomorphisms $\xymatrix{A\ar[r]^-\al&G\ar[r]^-\bt&B}$ with $G$ a $\Lb$-module is exact, then $\al$ and $\bt$ are both normal.
\vskip+2mm

{\bf  1.8.} \ Let $\xymatrix{G\ar[r]^-\al&Y\ar[r]^-\bt&Z}$ be a
sequence of  $\Lb$-semimodules and $\Lb$-homomorphisms with $G$ a
$\Lb$-module and $\bt\al=0$. Assume that the following is
satisfied: whenever $\bt(y_1)=\bt(y_2)$, $y_1,y_2\in Y$, one  has
$\al(g)+y_1=y_2$, $g\in G$. Then,  obviously, $\bt$ is a normal
$\Lb$-homomorphism and
    $\xymatrix{G\ar[r]^-\al&Y\ar[r]^-\bt&Z}$ is exact.
\vskip+2mm

{\bf  1.9.} \  {\bf  Lemma.} {\em Suppose given a commutative diagram of
$\Lb$-semi\-mo\-du\-les and $\Lb$-homomorphisms
$$\xymatrix{X\ar[r]^-\al\ar[d]_-f&Y\ar[r]^-\bt\ar[d]_\vf&Z\ar[d]_-\psi\\
X'\ar[r]^-{\al'}&Y'\ar[r]^-{\bt'}&Z'}$$
such that $f$ is surjective, $\vf$ is injective, and $\bt\al=0$. Assume that the bottom row
is exact and $\bt'$ is normal. Then the top row is exact and $\bt$ is normal.}

\begin{proof}
Suppose that $\bt(y_1)=\bt(y_2)$, $y_1,y_2\in Y$. Then $\bt'\vf(y_1)=\bt'\vf(y_2)$. Hence
$\kappa_1+\vf(y_1)=\kappa_2+\vf(y_2)$, $\kappa_1,\kappa_2\in\Ker(\bt')$. Since the bottom row
 is exact and $f$ is onto, there exist $x_1,x_2\in X$ such that $\kappa_1=\al'f(x_1)$ and
  $\kappa_2=\al'f(x_2)$. Then we get $\vf\al(x_1)+\vf(y_1)=\vf\al(x_2)+\vf(y_2)$. Whence, as
   $\vf$ is one-to-one,
$\al(x_1)+y_1=\al(x_2)+y_2$. Thus $\bt$ is normal. Now assume that $\bt(y)=0$, $y\in Y$.
Then $\bt'\vf(y)=0$. Hence $\al'f(x)=\vf(y)$ for some $x\in X$. This gives $\vf\al(x)=\vf(y)$. Whence $\al(x)=y$.
\end{proof}
\vskip+2mm

{\bf  1.10.} \ {\bf Definition ([4]).} We say that a sequence of $\Lb$-semimodules and
$\Lb$-homomorphisms
$$\xymatrix{X:\cdots\ar@<0.45ex>[r]\ar@<-0.45ex>[r]&X_{n+1}
\ar@<0.45ex>[r]^-{\pa_{n+1}^+}\ar@<-0.45ex>[r]_-{\pa_{n+1}^-}&
X_n\ar@<0.45ex>[r]^-{\pa_n^+}\ar@<-0.45ex>[r]_-{\pa_n^-}&
X_{n-1}\ar@<0.45ex>[r]\ar@<-0.45ex>[r]&\cdots},\quad n\in\bZ,$$
written $X=\{X_n,\pa_n^+,\pa_n^-\}$ for short, is a chain complex if
$$\pa_n^+\,\pa_{n+1}^++\pa_n^-\,\pa_{n+1}^-=\pa_n^+\,\pa_{n+1}^-+\pa_n^-\,
\pa_{n+1}^+$$
for each integer $n$. For every chain complex $X$ we define the $\Lb$-semi\-mo\-dule
$$Z_n(X)=\big\{x\in X_n|\pa_n^+(x)=\pa_n^-(x)\big\},$$
the $n$-cycles, and the $n$-th homology $\Lb$-semimodule
$$H_n(X)=Z_n(X)/\rho_n(X),$$
where $\rho_n(X)$ is a congruence on $Z_n(X)$ defined as follows:
\begin{align*}
x\rho_n(X)y\;\;\Leftrightarrow \;\; &x+\pa_{n+1}^+(u)+\pa_{n+1}^-(v)=y+\pa_{n+1}^+(v)+
\pa_{n+1}^-(u)\\
&\;\text{for some}\quad u,v \quad\text{in}\quad X_{n+1}.
\end{align*}
The $\Lb$-homomorphisms $\pa_n^+,\pa_n^-$ are called differentials of the chain complex $X$.
\vskip+4mm

A sequence $G=\{G_n,d_n^+,d_n^-\}$ of $\Lb$-modules and $\Lb$-homomorphisms is a chain
complex if and only if
$$\xymatrix{\cdots\ar[r]&G_n\ar[r]^-{d_n^+-d_n^-}&G_{n-1}\ar[r]&\cdots}$$
is an ordinary chain complex of $\Lb$-modules. Obviously, for any chain complex
$ G=\{G_n,d_n^+,d_n^-\}$ of $\Lb$-modules, $H_*(G)$ coincides with the usual homology
$H_*(\{G_n,d_n^+-d_n^-\})$.
\vskip+2mm

{\bf  1.11.} \  One can think of an ordinary chain complex
$$\xymatrix{\cdots\ar[r]&C_{n+1}\ar[r]^-{\pa_{n+1}}&C_n\ar[r]^-{\pa_n}&C_{n-1}\ar[r]
&\cdots}$$
of $\Lb$-semimodules as a chain complex in the sense of Definition 1.10; namely, we
identify $\{C_n,\pa_n\}$ with the chain complex
$$\xymatrix{\cdots\ar@<0.45ex>[r]\ar@<-0.45ex>[r]&C_{n+1}
\ar@<0.45ex>[r]^-{\pa_{n+1}}\ar@<-0.45ex>[r]_-0&
C_n\ar@<0.45ex>[r]^-{\pa_n}\ar@<-0.45ex>[r]_-0&
C_{n-1}\ar@<0.45ex>[r]\ar@<-0.45ex>[r]&\cdots}.$$
Defining $H_k(\{C_n,\pa_n\})$ to be $ H_k(\{C_n,\pa_n,0\})$, one has
$H_k(\{C_n,\pa_n\})=\Ker(\pa_k)/\pa_{k+1}(C_{k+1})$.
\vskip+2mm

{\bf  1.12.} \ {\bf  Definition ([4]).} Let $X=\{X_n,\pa_n^+,\pa_n^-\}$ and
$X'=\linebreak\{X'_n,{\pa_n^{\,'}}^+,{\pa_n^{\,'}}^-\}$ be chain complexes of
$\Lb$-semimodules. We say that a sequence $f=\{f_n\}$ of $\Lb$-homomorphisms $f_n:X_n\lra
X'_n$ is a $\pm$-morphism from $X$ to $X'$  if
$$f_{n-1}\pa_n^+={\pa_n^{\,'}}^+f_n\quad\text{and}\quad f_{n-1}\pa_n^-={\pa_n^{\,'}}^-f_n
\quad\text{for all}\quad n. $$
\vskip+2mm

{\bf  1.13.} \  If  $f=\{f_n\}:X\lra X'$ is a $\pm$-morphism of chain complexes, then
$f_n(Z_n(X))\sbs Z_n(X')$, and the map
$$H_n(f):H_n(X)\lra H_n(X'),\;\;H_n(f)(\cl(x))=\cl(f_n(x)),$$
is a homomorphism of $\Lb$-semimodules. Thus $H_n$ is a covariant additive functor from the
 category of chain complexes and their $\pm$-morphisms to the category of $\Lb$-semimodules.
\vskip+4mm

An important example of a $\pm$-morphism appears in a natural way: a map of presimplicial
 $\Lb$-semimodules $f:S\lra S'$ induces a $\pm$-mor\-phism $\ul{f}=f:\ul{S}\lra\ul{S}'$,
 where $\ul{S}$ and $\ul{S}'$ are the standard nonnegative chain complexes associated to
 $S$ and $S'$, respectively (see [4]).
\vskip+4mm

{\bf  1.14.} \  {\bf Definition (cf. [5]).} Let $X=\{X_n,\pa_n^+,\pa_n^-\}$ and
$ X'=\{X'_n,{\pa_n^{\,'}}^+,{\pa_n^{\,'}}^-\}$ be chain complexes of $\Lb$-semimodules. We
say that a sequence $f=\{f_n\}$ of $\Lb$-homomorphisms $f_n:X_n\lra X'_n$ is a morphism from
$X$ to $X'$ if
$${\pa_n^{\,'}}^+f_n+f_{n-1}\pa_n^-={\pa_n^{\,'}}^-f_n+f_{n-1}\pa_n^+\quad\text{for all}
\quad n.$$
\vskip+4mm

Note that any $\pm$-morphism between chain complexes of $\Lb$-semi\-mo\-du\-les is a morphism.
\vskip+2mm

{\bf  1.15.} \ {\bf Definition.} A sequence $\xymatrix{E:A\ar@{ >->}[r]^-{\vkp}&B\ar@{->>}
[r]^-{\sg}&C}$ of chain complexes and their morphisms is called a Schreier short exact
sequence of chain complexes if each
$\xymatrix{E_n:A_n\ar@{ >-{>}}[r]^-{\vkp_n}&B_n\ar@{ ->>}[r]^-{\sg_n}&C_n}$ is a Schreier
extension of $\Lb$-semimodules.
\vskip+2mm

{\bf  1.16.} \  In general, a morphism $f=\{f_n\}:X\lra X'$ of chain complexes, unlike
$\pm$-morphisms, does not induce a $\Lb$-homomorphism from $H_n(X)$ to $H_n(X')$. However, if
$X'$ is a chain complex of cancellative $\Lb$-semimodules, or $X$ is an ordinary chain complex of $\Lb$-semi\-mo\-dules, i.e., $\pa_n^-=0$ for all $n$ (see 1.11), then one can easily check that $f_n(Z_n(X))\sbs Z_n(X')$ and the map $H_n(f):H_n(X)\lra H_n(X')$,
$H_n(f)(\cl(x))=\cl(f_n(x))$, is well-defined and is a $\Lb$-homo\-mor\-phism for all $n$.
Besides, we have
\vskip+2mm

{\bf Proposition.} {\em  Let $\xymatrix{E:A\ar@{ >-{>}}[r]^-\vkp&B\ar@{ ->>}[r]^-\sg&C}$ be
a Schreier short exact sequence of chain complexes and their morphisms. If $A$ is a chain
complex of cancellative $\Lb$-semimodules, then
$\vkp_n(Z_n(A))\sbs Z_n(B)$ and the map $H_n(\vkp):H_n(A)\lra H_n(B)$, $H_n(\vkp)(\cl(a))=
\cl(\vkp_n(a))$, is well-defined and is therefore a $\Lb$-homo\-mor\-phism for all $n$.}
\vskip+2mm

\begin{proof}
Let $d_n^+,d_n^-$ and $\pa_n^+,\pa_n^-$ denote the $n$-th differentials of $A$ and $B$,
respectively. Suppose $a\in Z_n(A)$, i.e., $d_n^+(a)=d_n^-(a)$. Since $\vkp$ is a morphism,
$\vkp_{n-1}d_n^+(a)+\pa_n^-\vkp_n(a)=\vkp_{n-1}d_n^-(a)+\pa_n^+\vkp_n(a)$. Whence,by 1.3,
 $\pa_n^-\vkp_n(a)=\pa_n^+\vkp_n(a)$. That is, $\vkp_n(a)\in Z_n(B)$. Now assume that
$a_1,a_2\in Z_n(A)$ and $a_1\rho_n(A)a_2$. Hence
$$a_1+d_{n+1}^+(p)+d_{n+1}^-(q)=a_2+d_{n+1}^+(q)+d_{n+1}^-(p),\quad p,q\in A_{n+1}.$$
On the other hand,
\begin{align*}
\vkp_nd_{n+1}^+(p)+\pa_{n+1}^-\vkp_{n+1}(p)&=\vkp_nd_{n+1}^-(p)+
\pa_{n+1}^+\vkp_{n+1}(p),\\
\vkp_nd_{n+1}^+(q)+\pa_{n+1}^-\vkp_{n+1}(q)&=\vkp_nd_{n+1}^-(q)+
\pa_{n+1}^+\vkp_{n+1}(q).
\end{align*}
These last three equalities imply
\begin{gather*}
\vkp_n(d_{n+1}^+(p)+d_{n+1}^-(q))+\pa_{n+1}^+\vkp_{n+1}(p)+\pa_{n+1}^-\vkp_{n+1}(q)+
\vkp_n(a_1)=\\
=\vkp_n(d_{n+1}^+(p)+d_{n+1}^-(q))+\pa_{n+1}^+\vkp_{n+1}(q)+\pa_{n+1}^-\vkp_{n+1}(p)+
\vkp_n(a_2).
 \end{gather*}
Whence, by 1.3, $\pa_{n+1}^+\vkp_{n+1}(p)+\pa_{n+1}^-\vkp_{n+1}(q)+\vkp_n(a_1)=
\pa_{n+1}^+\vkp_{n+1}(q)+\pa_{n+1}^-\vkp_{n+1}(p)+\vkp_n(a_2)$. That is,
$\vkp_n(a_1)\rho_n(B)\vkp_n(a_2)$. Thus $H_n(\vkp)$ is well-defined.
\end{proof}
\vskip+4mm

 Definition 1.10  naturally leads us to new homology and cohomology monoids of monoids
 (in particular, groups) with coefficients in semimodules. The calculation of them for cyclic
  groups is an example of the effective use of morphisms which are not $\pm$-morphisms [6].
\vskip+4mm

{\bf  1.17.} \ If $X=\{X_n,\pa_n^+,\pa_n^-\}$ is a chain complex of $\Lb$-semimodules, then
$K(X)=\{K(X_n),K(\pa_n^+)-K(\pa_n^-)\}$ is an ordinary chain complex of $K(\Lb)$-modules
(i.e., $\Lb$-modules). When each $X_n$ is cancellative, then the converse is also true. Then,
 for any chain complex $X=\{X_n,\pa_n^+,\pa_n^-\}$ of $\Lb$-semimodules, one has the
 $\Lb$-homomorphisms $H_n(k_X):H_n(X)\lra H_n(K(X))$, $H_n(k_X)(\cl(x))=\cl(k_{X_n}(x))=
 \cl[x,0]$, induced by the canonical morphism $k_X=\{k_{X_n}\}:X\lra K(X)$ which is in fact
 a $\pm$-morphism from $X$ to $\{K(X_n),K(\pa_n^+),K(\pa_n^-)\}$. When $X$ is a chain complex
  of cancellative
$\Lb$-semimodules, then $H_n(k_X)$ is injective and therefore $H_n(X)$ is a cancellative
$\Lb$-semimodule. Further,  if  $f=\{f_n\}:X\lra X'$ is a morphism of chain complexes, then
 $K(f)=\{K(f_n):K(X_n)\lra K(X'_n)\}$ is an usual chain map from $K(X)$ to  $K(X')$. When
 $X'$ is a chain complex of cancellative $\Lb$-semimoduls, then the converse is also valid.
\vskip+2mm

{\bf 1.18.} \ If $\xymatrix{E:A\ar@{ >-{>}}[r]^-\vkp&B\ar@{ ->>}[r]^-\sg&C}$ is a Schreier
short exact sequence of chain complexes, then, by 1.5,
$\xymatrix{K(E):K(A)\ar@{ >-{>}}[r]^-{K(\vkp)}&K(B)\ar@{ ->>}[r]^-{K(\sg)}&K(C)}$ is a short
 exact sequence of ordinary chain complexes.

\section{Main results}

{\bf 2.1.} \  {\bf Proposition.} {\em Suppose given a Schreier short exact sequence
$$\xymatrix{A\ar@{ >-{>}}[r]^-\vkp&B\ar@{ ->>}[r]^-\sg&C}$$
of chain complexes and their morphisms such that each $A_n$ is cancellative and each
differential $\pa_n^-$ of $B$ preserves representatives. Assume that one of the following
conditions holds:

\ {\rm (i)} \ $\sg$ is  a $\pm$-morphism.

{\rm (ii)} \ $C$ is a chain complex of cancellative $\Lb$-semimodules.

\noindent Then there are $\Lb$-homomorphisms $\pa_n(E):H_n(C)\lra H_{n-1}(A)$, called
 connecting homomorphisms, such that each diagram

$$\xymatrix{H_n(C)\ar[rr]^-{\pa_n(E)}\ar[dd]_-{H_n(k_C)}&&H_{n-1}(A)
\ar[dd]^-{H_{n-1}(k_A)}\\
&&\\
H_n(K(C))\ar[rr]^{\pa_n(K(E))}&&H_{n-1}(K(A)),}$$
where $\pa_n(K(E))$ is the usual connecting homomorphism induced by $K(E)$ $($see $1.18)$,
is commutative. Furthermore, $H_n(\vkp)$ and $H_n(\sg)$ are defined for all $n$, and the
long sequence of homology semimodules
$$\!\xymatrix{\!\cdots\!\to\! H_n(A)\ar[r]^-{H_n(\vkp)}&H_n(B)\ar[r]^-{H_n(\sg)}&
H_n(C)\ar[r]^-{\pa_n(E)}&H_{n\!-\!1}(A)\ar[r]^-{H_{n-1}(\vkp)}
&H_{n\!-\!1}(B)\!\to\!\cdots}$$
is an ordinary chain complex of $\Lb$-semimodules.}
\vskip+2mm

\begin{proof}
Let $d_n^+,d_n^-$ and $\dl_n^+,\dl_n^-$ denote the $n$-th differentials of $A$ and $C$,
respectively.  Take any $c\in Z_n(C)$. There is a representative $u_c$ of
$\xymatrix{E_n:A_n\ar@{ >-{>}}[r]^-{\vkp_n}&B_n\ar@{ ->>}[r]^-{\sg_n}&C_n}$ with $\sg_n(u_c)
=c$. When (i) holds, one can write $\sg_{n-1}\pa_n^+(u_c)=\dl_n^+\sg_n(u_c)=
\dl_n^+(c)=\dl_n^-(c)=\dl_n^-\sg_n(c)=\sg_{n-1}\pa_n^-(u_c)$. If (ii)  holds, then the
equality $\sg_{n-1}\pa_n^+(u_c)+\dl_n^-\sg_n(u_c)=\sg_{n-1}\pa_n^-(u_c)+\dl_n^+\sg_n(u_c)$
implies  $\sg_{n-1}\pa_n^+(u_c)=\sg_{n-1}\pa_n^-(u_c)$. Consequently, $\pa_n^+(u_c)=\vkp_{n-1}
(a)+\pa_n^-(u_c)$, $a\in A_{n-1}$, in both cases (see 1.1). Whence $[\pa_n^+(u_c),\pa_n^-
(u_c)]=K(\vkp_{n-1})([a,0])$. On the other hand, $[\pa_n^+(u_c),\pa_n^-(u_c)]\linebreak=
(K(\pa_n^+)-K(\pa_n^-))([u_c,0])$ and $K(\sg_n)([u_c,0])=[c,0]\in\Ker(K(\dl_n^+)-K(\dl_n^-))$. Therefore, by construction of  $\pa_n(K(E))$, one concludes that $[a,0]\in\Ker(K(d_{n-1}^+)-K(d_{n-1}^-))$  and $\pa_n(K(E))(\cl([c,0]))=\cl([a,0])$. As $A_{n-2}$ is cancellative, the former gives $a\in Z_{n-1}(A)$. And we set
$$\pa_n(E)(\cl(c))=\cl(a)\in H_{n-1}(A).$$
Clearly, since $\pa_n(K(E))H_n(k_C)(\cl(c))=H_{n-1}(k_A)(\cl(a))$ and \linebreak$H_{n-1}(k_A)$ is injective (see 1.17), we may write
$$\pa_n(E)(\cl(c))=H_{n-1}(k_A)^{-1}\big(\pa_n(K(E))H_n(k_C)(\cl(c))\big).$$
Hence, $\pa_n(E)$ is well-defined and is a $\Lb$-homomorphism, and \linebreak $H_{n-1}(k_A)
\pa_n(E)=\pa_n(K(E))H_n(k_C)$.

It follows from  1.13  and 1.16  that
$H_n(\sg)$ and $H_n(\vkp)$ are defined. Obviously $H_n(\sg)H_n(\vkp)=0$. Using the usual
long exact homology sequence for $K(E)$, one has \begin{gather*}
H_{n-1}(k_A)\pa_n(E)H_n(\sg)=\pa_n(K(E))H_n(k_C)
H_n(\sg)=\\
=\pa_n(K(E))H_n(K(\sg))H_n(k_B)=0\cdot H_n(k_B)=0.
\end{gather*}
Hence  $\pa_n(E)H_n(\sg)=0$ since $H_{n-1}(k_A)$ is one-to-one. By definition of $\pa_n(E)$,
 $\pa_n(E)(\cl(c))=\cl(a)$, $a\in Z_{n-1}(A)$, and $a$ satisfies the equality $\pa_n^+(u_c)=
 \vkp_{n-1}(a)+\pa_n^-(u_c)$ for some representative $u_c$ of $E_n$ with $\sg_n(u_c)=c$.
 Consequently, $H_{n-1}(\vkp)\pa_n(E)(\cl(c))=H_{n-1}(\vkp)(\cl(a))=\cl(\vkp_{n-1}(a))=0$.
  Thus $H_{n-1}(\vkp)\pa_n(E)=0$.
\end{proof}
\vskip+2mm

We see that for any Schreier short exact sequence of chain complexes
$\xymatrix{E:A\ar@{ >-{>}}[r]^-\vkp&B\ar@{ ->>}[r]^-\sg&C}$ satisfying the hypotheses of
Proposition 2.1, one has the commutative diagram
{\small $$\xymatrix{\cdots\to H_n(A)\ar[r]^-{H_n(\vkp)}\ar[dd]_-{H_n(k_A)}&
H_n(B)\ar[r]^-{H_n(\sg)}\ar[dd]_-{H_n(k_B)}&H_n(C)\ar[r]^-{\pa_n(E)}
\ar[dd]_-{H_n(k_C)}&H_{n-1}(A)\rightarrow\ar[dd]_-{H_{n-1}(k_A)}\cdots\\
&&&\\
\!\cdots\!\to\! H_n(K(A))\;\ar[r]^-{H_n(K(\vkp))}&\;H_n(K(B))\;\ar[r]^-{H_n(K(\sg))}
&\;H_n(K(C))\;\ar[r]^-{\pa_n(K(E))}&\;H_{n\!-\!1}(K(A))\!\rightarrow\!\cdots \!}$$}
induced by the canonical map $k_E=(k_A,k_B,k_C):E\lra K(E)$. In fact $\pa_n(E)$ is natural in
 the following sense.
\vskip+2mm

{\bf 2.2.} \  {\bf Proposition.} {\em Suppose that
$$\xymatrix{E:\!\!\!&\!\!\!A\;\ar@{ >->}[r]^-\vkp\ar[d]_-f&B\ar@{ ->>}[r]^\sg\ar[d]_-
g&C\ar[d]_-h\\
E':\!\!\!&\!\!\!A'\;\ar@{ >->}[r]^-{\vkp'}&B'\ar@{ ->>}[r]^-{\sg'}&C'}$$
is a commutative diagram of chain complexes and their morphisms such that $E$ and $E'$ are
Schreier short exact sequences satisfying the hypotheses of Proposition $2.1$. Suppose
further that $H_n(g)$ and $H_n(h)$ are defined for all $n$ $($see $1.13$ and $1.16)$. Then
the diagram
{\small $$\xymatrix{\cdots\to H_n(A)\ar[r]^-{H_n(\vkp)}\ar[dd]_-{H_n(f)}&
H_n(B)\ar[r]^-{H_n(\sg)}\ar[dd]_-{H_n(g)}&H_n(C)\ar[r]^-{\pa_n(E)}
\ar[dd]_-{H_n(h)}&H_{n-1}(A)\rightarrow\ar[dd]_-{H_{n-1}(f)}\cdots\\
&&&\\
\cdots\to H_n(A')\;\ar[r]^-{H_n(\vkp')}&\;H_n(B')\;\ar[r]^-{H_n(\sg')}
&\;H_n(C')\;\ar[r]^-{\pa_n(E')}&\;H_{n-1}(A')\rightarrow\cdots }$$}
is commutative.}
\vskip+2mm

\begin{proof}
By  Proposition 2.1 and the naturality of $\pa_n(K(E))$,
\begin{align*}
&H_{n-1}(k_{A'})\pa_n(E')H_n(h)=\pa_n(K(E'))H_n(k_{C'})H_n(h)\\
&=\!\pa_n(K(E'))H_n(K(h))H_n(k_C)\!=\!H_{n-1}(K(f))\pa_n(K(E))H_n(k_C)\\
&=H_{n-1}(K(f))H_{n-1}(k_A)\pa_n(E)=H_{n-1}(k_{A'})H_{n-1}(f)\pa_n(E).
\end{align*}
Therefore, by the injectivity of $H_{n-1}(k_{A'})$ (see 1.17), $\pa_n(E')H_n(h)=H_{n-1}(f)
\pa_n(E)$.
\end{proof}

Before we state our main results, we note the following. Let
\linebreak $\xymatrix{A\ar[r]^- \al&B\ar[r]^-\bt&C}$ be an exact
sequence of $\Lb$-semi\-mo\-dules and $\Lb$-ho\-mo\-mor\-phisms.
If $\bt=0$ then $\al$ is onto. But, unlike the situation for
modules, one may have
 $\al=0$ and yet not have $\bt$ one-to-one. However, we have:
\vskip+2mm

{\bf 2.3.} \  Suppose given an exact sequence of $\Lb$-semimodules
and $\Lb$-ho\-mo\-morphisms
$\xymatrix{A\ar[r]^-\al&B\ar[r]^-\bt&C}$ with $\bt$ a normal
$\Lb$-homomorphism (see 1.6). If $\al=0$, then $\bt$ is
one-to-one. \vskip+2mm

This together with 1.7 motivates the following two theorems.
\vskip+2mm

{\bf 2.4.} \  {\bf Theorem.} {\em  Let
{\small$$\xymatrix{\;\;\;\;\;\;A:\cdots\ar@{
>->}[dd]^-\vkp\ar@<0.45ex>[r]\ar@<-0.45ex>[r] &A_{n+1}
\ar@<0.45ex>[rr]^-{d_{n+1}^+}\ar@<-0.45ex>[rr]_-{d_{n+1}^-}\ar@{
>->}[dd]^-{\vkp_{n+1}}&&
A_n\ar@<0.45ex>[rr]^-{d_n^+}\ar@<-0.45ex>[rr]_-{d_n^-}\ar@{
>->}[dd]^-{\vkp_n}&&A_{n-1}
\ar@<0.45ex>[r]\ar@<-0.45ex>[r]\ar@{ >->}[dd]^-{\vkp_{n-1}}&\cdots\\
&&&\\
E:\;\;\;\;\;\;B:\;\;\;\;\cdots\ar@{
->>}[dd]^-\sg\ar@<0.45ex>[r]\ar@<-0.45ex>[r]&B_{n+1}
\ar@<0.45ex>[rr]^-{\pa_{n+1}^+}\ar@<-0.45ex>[rr]_-{\pa_{n+1}^-}\ar@{
->>}[dd]^-{\sg_{n+1}}
&&B_n\ar@<0.45ex>[rr]^-{\pa_n^+}\ar@<-0.45ex>[rr]_-{\pa_n^-}\ar@{
->>}[dd]^-{\sg_n}&&B_{n-1}
\ar@<0.45ex>[r]\ar@<-0.45ex>[r]\ar@{ ->>}[dd]^-{\sg_{n-1}}&\cdots\\
&&&\\
\;\;\;\;\;\;G:\cdots\ar[r]&G_{n+1}\ar[rr]^-{\dl_{n+1}}&&G_n
\ar[rr]^-{\dl_n}&&G_{n-1}\ar[r]&\cdots}$$}
be a Schreier short exact sequence of chain complexes and their morphisms $($see $1.11)$ such
that each $A_n$ is a cancellative $\Lb$-semimodule, every differential $\pa_n^-$ preserves
 representatives, and each $G_n$ is a $\Lb$-mo\-du\-le. Then the long homology sequence
\begin{gather}
\xymatrix{\cdots\!\to\! H_{n+1}(G)\ar[r]^-{\pa_{n+1}(E)}
&H_n(A)\ar[r]^-{H_n(\vkp)}&H_n(B)\ar[r]^-{H_n(\sg)}&}\tag{$2.4.1$}\\
\xymatrix{{} \ar[r]^-{H_n(\sg)}&H_n(G)\ar[r]^-{\pa_n(E)}&H_{n-1}(A)\ar[r]&\cdots}\notag\
\end{gather}
is exact at $H_n(A)$ and at $H_n(B)$, $H_n(\sg)(H_n(B))\sbs\Ker(\pa_n(E))$, and
$H_n(\vkp)$ is normal. Furthermore, $(2.4.1)$ is exact at $H_n(G)$ if and only if
$H_n(\sg)(H_n(B))=H_n(K(\sg))(H_n(K(B)))$.}
\vskip+2mm

\begin{proof}
By Proposition 2.1, Sequence (2.4.1) is an ordinary chain complex. Then, the commutative
diagram
$$\xymatrix{H_{n+1}(G)\ar[rr]^{\pa_{n+1}(E)}\ar@{=}[dd]&&H_n(A)
\ar[rr]^-{H_n(\vkp)}
\ar[dd]^-{H_n(k_A)}&&H_n(B)\ar[dd]^-{H_n(k_B)}\\
&&\\
H_{n+1}(G)\ar[rr]^-{\pa_{n+1}(K(E))}&&H_n(K(A))\ar[rr]^-{H_n(K(\vkp))}
&&H_n(K(B))}$$
satisfies the hypotheses of Lemma 1.9 (see 1.17). Hence (2.4.1) is exact at $H_n(A)$ and
 $H_n(\vkp)$ is normal. We next show that $\Ker(H_n(\sg))\sbs H_n(\vkp)(H_n(A))$. Let
 $b\in Z_n(B)$, i.e., $\pa_n^+(b)=\pa_n^-(b)$. Assume that $H_n(\sg)(\cl(b))=0$, i.e.,
 $\cl(\sg_n(b))=0$. Then $\sg_n(b)=\dl_{n+1}(g)$ for some \linebreak $g\in G_{n+1}$. Choose
 a representative $u=u_{-g}$ of
$$\xymatrix{E_{n+1}:A_{n+1}\ar@{ >->}[r]^-{\vkp_{n+1}}&B_{n+1}\ar@{ ->>}[r]^-{\sg_{n+1}}&
G_{n+1}}$$ with $\sg_{n+1}(u)=-g$. Since $\sg$ is a morphism,
$\sg_n\pa_{n+1}^+(u)=\sg_n\pa_{n+1}^-(u)+\dl_{n+1}\sg_{n+1}(u)$.
Whence $\sg_n\pa_{n+1}^+(u)=\sg_n\pa_{n+1}^-(u)-\sg_n(b)$, i.e.,
$\sg_n(\pa_{n+1}^+(u)+b)=\sg_n\pa_{n+1}^-(u)$. Then, as
$\pa_{n+1}^-$ preserves representatives, if follows that
\begin{equation}
b+\pa_{n+1}^+(u)=\vkp_n(a)+\pa_{n+1}^-(u),\quad a\in A_n.\tag{$*$}
\end{equation}
This  with the fact that $\vkp$ is a morphism of chain complex gives
\begin{align*}
\vkp_{n-1}d_n^+(a)+\pa_n^-(b)+(\pa_n^-\pa_{n+1}^++\pa_n^+\pa_{n+1}^-)(u)=\\
=\vkp_{n-1}d_n^-(a)+\pa_n^+(b)+(\pa_n^+\pa_{n+1}^++\pa_n^-\pa_{n+1}^-)(u).
\end{align*}
But, by 1.4, $B_{n-1}$ is cancellative. Therefore $\vkp_{n-1}d_n^+(a)=\vkp_{n-1}d_n^-(a)$
(see 1.10). Hence $d_n^+(a)=d_n^-(a)$, i.e., $a\in Z_n(A)$. Then, by $(*)$, one can write
$H_n(\vkp)(\cl(a))=\cl(\vkp_n(a))=\cl(b)$. Thus (2.4.1) is exact at $H_n(B)$. Finally, the
 commutative diagram
$$\xymatrix{H_n(B)\ar[rr]^-{H_n(\sg)}\ar[dd]_-{H_n(k_B)}&&H_n(G)\ar[rr]^-{\pa_n(E)}
\ar@{=}[dd]&&H_{n-1}(A)\ar[dd]^-{H_{n-1}(k_A)}\\
&&\\
H_n(K(B))\ar[rr]^{H_n(K(B))}&&H_n(G)\ar[rr]^-{\pa_n(K(E))}&&H_{n-1}(K(A))}$$
shows that if $H_n(\sg)(H_n(B))=H_n(K(\sg))(H_n(K(B)))$, then (2.4.1) is exact at $H_n(G)$.
 The converse is also true since $H_{n-1}(k_A)$ is injective (see 1.17).
\end{proof}
\vskip+2mm

{\bf 2.5. Theorem.} {\em Suppose given a Schreier short exact sequence
{\small$$\xymatrix{\;\;\;\;\;\;\;\;\;G:\;\;\cdots\ar[r]\ar@{ >->}[dd]^-{\vkp}
&G_{n+1}\ar[rr]^-{d_{n+1}}\ar@{ >->}[dd]^-{\vkp_{n+1}}
&&G_n\ar[rr]^-{d_n}\ar@{>->}[dd]^-{\vkp_n}
&&G_{n-1}\ar[r]\ar@{ >->}[dd]^-{\vkp_{n-1}}&\cdots\\
&&&\\
E:\;\;\;\;\;B:\;\;\;\cdots\ar@{
->>}[dd]^-\sg\ar@<0.45ex>[r]\ar@<-0.45ex>[r]&B_{n+1}
\ar@<0.45ex>[rr]^-{\pa_{n+1}^+}\ar@<-0.45ex>[rr]_-{\pa_{n+1}^-}\ar@{
->>}[dd]^-{\sg_{n+1}}
&&B_n\ar@<0.45ex>[rr]^-{\pa_n^+}\ar@<-0.45ex>[rr]_-{\pa_n^-}\ar@{
->>}[dd]^-{\sg_n}&&B_{n-1}
\ar@<0.45ex>[r]\ar@<-0.45ex>[r]\ar@{ ->>}[dd]^-{\sg_{n-1}}&\cdots\\
&&&\\
\;\;\;\;\;\;\;\;\;C:\;\;\cdots\ar@<0.45ex>[r]\ar@<-0.45ex>[r]
&C_{n+1}\ar@<0.45ex>[rr]^-{\dl_{n+1}^+}\ar@<-0.45ex>[rr]_-{\dl_{n+1}^-}
&&C_n\ar@<0.45ex>[rr]^-{\dl_n^+}\ar@<-0.45ex>[rr]_-{\dl_n^-}
&&C_{n-1}\ar@<0.45ex>[r]\ar@<-0.45ex>[r]&\cdots}$$}
of chain complexes and their morphisms $($see $1.11)$ such that each $G_n$ is a $\Lb$-module
$($therefore each differential $\pa_n^-$ obviously preserves representatives $($see $1.2))$.
 Assume that one of the following conditions holds:

\ {\rm (i)} \ $\sg$ is a $\pm$-morphism.

{\rm (ii)} \ $C$ is a chain complex of cancellative $\Lb$-semimodules.

\noindent Then the long homology sequence
\begin{gather}
\xymatrix{\cdots\!\to\! H_n(G)\ar[r]^-{H_n(\vkp)}
&H_n(B)\ar[r]^-{H_n(\sg)}&H_n(C)\ar[r]^-{\pa_n(E)}&}\tag{$2.5.1$}\\
\;\;\;\xymatrix{{}\ar[r]^-{\pa_n(E)}&H_{n-1}(G)\ar[rr]^-{H_{n-1}(\vkp)}
&&H_{n-1}(B)\ar[r]&
\cdots}\notag
\end{gather}
$($here $\pa_n(E)$ evidently coincides with $\pa_n(K(E))H_n(k_C))$ is exact at \linebreak
$H_n(B)$ and at $H_n(C)$, $\pa_n(E)(H_n(C))\sbs \Ker(H_{n-1}(\vkp))$, and $H_n(\sg)$ is a
normal $\Lb$-homomorphism. Furthermore, if \ $\pa_n(K(E))H_n(k_C)(H_n(C))=\pa_n(K(E))
(H_n(K(C)))$ then $(2.5.1)$ is exact at $H_{n-1}(G)$. When {\rm (ii)} holds, the converse
is also valid.}
\vskip+2mm

\begin{proof}
By Proposition 2.1, Sequence (2.5.1) is an ordinary chain complex. Assume that (ii) holds.
Then, by 1.4, the commutative diagram
$$\xymatrix{H_n(G)\ar[rr]^{H_n(\vkp)}\ar@{=}[dd]&&H_n(B)\ar[rr]^-{H_n(\sg)}
\ar[dd]^-{H_n(k_B)}&&H_n(C)\ar[dd]^-{H_n(k_C)}\\
&&\\
H_n(G)\ar[rr]^-{H_n(K(\vkp))}&&H_n(K(B))\ar[rr]^-{H_n(K(\sg))}&&H_n(K(C))}$$
satisfies the hypotheses of Lemma 1.9 (see 1.17).   Therefore (2.5.1) is exact at $H_n(B)$
and $H_n(\sg)$ is normal. When (i) holds, we prove the same as follows. Let
$b_1,b_2\in Z_n(B)$, i.e., $\pa_n^+(b_1)=\pa_n^-(b_1)$ and $\pa_n^+(b_2)=\pa_n^-(b_2)$, and
 let $H_n(\sg)(\cl(b_1))=H_n(\sg)(\cl(b_2))$, i.e., $\cl(\sg_n(b_1))\!\!=\!\cl(\sg_n(b_2))$.
  Then $\sg_n(b_1)+\dl_{n+1}^+(p)+\dl_{n+1}^-(q)=\sg_n(b_2)+\dl_{n+1}^+(q)+\dl_{n+1}^-(p)$
   for some $p,q\in C_{n+1}$. Take  $x,y\in B_{n+1}$ with $\sg_{n+1}(x)=p$ and
   $\sg_{n+1}(y)=q$, and write
$\sg_n(b_1)+\dl_{n+1}^+\sg_{n+1}(x)+\dl_{n+1}^-\sg_{n+1}(y)=\sg_n(b_2)+\dl_{n+1}^+
\sg_{n+1}(y)+\dl_{n+1}^-\sg_{n+1}(x)$. This, as $\sg$ is a $\pm$-morphism, implies
$$\sg_n\big(b_1+\pa_{n+1}^+(x)+\pa_{n+1}^-(y)\big)=\sg_n\big(b_2+\pa_{n+1}^+(y)+
\pa_{n+1}^-(x)\big).$$
Therefore, by 1.2,
\begin{equation}
b_1+\pa_{n+1}^+(x)+\pa_{n+1}^-(y)=\vkp_n(g)+b_2+\pa_{n+1}^+(y)+
\pa_{n+1}^-(x),\quad g\in G_n.\tag{$**$}
\end{equation}
Whence
\begin{gather*}
\pa_n^+(b_1)+\pa_n^+\pa_{n+1}^+(x)+\pa_n^+\pa_{n+1}^-(y)=\\
=\pa_n^+\vkp_n(g)+\pa_n^+b_2+\pa_n^+\pa_{n+1}^+(y)+\pa_n^+\pa_{n+1}^-(x)
\end{gather*}
and
\begin{gather*}
\pa_n^-(b_1)+\pa_n^-\pa_{n+1}^+(x)+\pa_n^-\pa_{n+1}^-(y)=\\
=\pa_n^-\vkp_n(g)+\pa_n^-(b_2)+\pa_n^-\pa_{n+1}^+(y)+\pa_n^-\pa_{n+1}^-(x).
\end{gather*}
The last two equalities give
\begin{align*}
\pa_n^+&\vkp_n(g)+\pa_n^+(b_2)+\pa_n^+\pa_{n+1}^+(y)+\pa_n^+\pa_{n+1}^-(x)+\\
&+\pa_n^-(b_1)+\pa_n^-\pa_n^+(x)+\pa_n^-\pa_{n+1}^-(y)=\\
=&\pa_n^-\vkp_n(g)+\pa_n^-(b_2)+\pa_n^-\pa_{n+1}^+(y)+\pa_n^-\pa_{n+1}^-(x)+\\
&+\pa_n^+(b_1)+\pa_n^+\pa_n^+(x)+\pa_n^+\pa_{n+1}^-(y).
\end{align*}
But $\pa_n^+(b_1)=\pa_n^-(b_1)$, $\pa_n^+(b_2)=\pa_n^-(b_2)$,
$\pa_n^+\vkp_n(g)=\vkp_{n-1}d_n(g)+\pa_n^-\vkp_n(g)$ and
$\pa_n^+\pa_{n+1}^++\pa_n^-\pa_{n+1}^-=\pa_n^+\pa_{n+1}^-+\pa_n^-\pa_{n+1}^+$.
Consequently, we have
$$\vkp_{n-1}d_n(g)+w=w,\quad w\in B_{n-1}.$$
Whence, by 1.2, $d_n(g)=0$. That is, $g\in\Ker(d_n)$. Then, by $(**)$, one can write
$$\cl(b_1)=H_n(\vkp)(\cl(g))+\cl(b_2).$$
Thus, by 1.8, we conclude that (2.5.1) is exact at $H_n(B)$ and $H_n(\sg)$ is normal.

We next show that $\Ker(\pa_n(E))\sbs H_n(\sg)(H_n(B))$. Let $c\in Z_n(C)$. Take any
$b\in B_n$ with $\sg_n(b)=c$. By definition of $\pa_n(E)$, $\pa_n(E)(\cl(c))=\cl(g)$,
$g\in Z_{n-1}(G)$, and $g$ satisfies the equality $\pa_n^+(b)=\vkp_{n-1}(g)+\pa_n^-(b)$
(see 1.2). Assume that $\pa_n(E)(\cl(c))=0$, i.e., $\cl(g)=0$. Then $g=d_n(h)$, $h\in G_n$.
As $\vkp$ is a morphism, we can write $\pa_n^+(b-\vkp_n(h))=\pa_n^+(b)-\pa_n^+\vkp_n(h)=
\vkp_{n-1}(g)+\pa_n^-(b)-
(\vkp_{n-1}d_n(h)+\pa_n^-\vkp_n(h))=\pa_n^-(b)-\pa_n^-\vkp_n(h)=\pa_n^-(b-\vkp_n(h))$. Hence
 $b-\vkp_n(h)\in Z_n(B)$. Clearly, $H_n(\sg)(\cl(b-\vkp_n(h)))=\cl(c)$. Thus (2.5.1) is exact
 at $H_n(C)$.

Finally, the commutative diagram
$$\xymatrix{H_n(C)\ar[rr]^-{\pa_n(E)}\ar[dd]_-{H_n(k_C)}&&H_{n-1}(G)
\ar[rr]^-{H_{n-1}(\vkp)}\ar@{=}[dd]&&H_{n-1}(B)\ar[dd]^-{H_{n-1}(k_B)}\\
&&\\
H_n(K(C))\ar[rr]^{\pa_n(K(E))}&&H_{n-1}(G)\ar[rr]^-{H_{n-1}(K(\vkp))}&&H_{n-1}(K(B))}$$
shows that if $\pa_n(K(E))H_n(k_C)(H_n(C))=\pa_n(K(E))(H_n(K(C)))$, then (2.5.1) is exact at
 $H_{n-1}(G)$. When (ii) holds, the converse is also true since $H_{n-1}(k_B)$ is injective
 (see 1.4 and 1.17).
\end{proof}
\vskip+2mm

{\bf 2.6.} {\bf Remark}. For a Schreier short exact sequence $E:$ \linebreak
$\xymatrix{G\ar@{ >->}[r]^-\vkp&B\ar@{ ->>}[r]^-\sg&C}$ of chain complexes, where $G$ is an
 ordinary chain complex of $\Lb$-modules, one always has the connecting homomorphism $\pa_n
 (K(E))H_n(k_C)$. But in general $H_{n-1}(\vkp)\pa_n(K(E))H_n(k_C)\neq 0$. Moreover, if
  neither (i) nor (ii) holds, Theorem 2.5 need not hold even in the case when $H_{n-1}(\vkp)
  \pa_n(K(E))H_n(k_C)=0$ and $H_n(\sg)$ is defined for all $n$ (see 1.16). Indeed, consider
  the following diagram
$$\xymatrix{\;G:\ar@{ >->}[d]_-\vkp&\cdots\ar[r]&0\ar[r]\ar@{ >->}[d]&0\ar[r]\ar@{ >->}[d]&
0\ar[r]\ar@{ >->}[d]&0\ar[r]\ar@{ >->}[d]&0\ar[r]\ar@{ >->}[d]&\cdots\\
\;B:\ar@{ >->}[d]_-\sg&\cdots\ar[r]&0\ar[r]\ar@{ >->}[d]&M\ar[r]^-1\ar@{ >->}[d]_-1
&M\ar[r]^-0\ar@{ >->}[d]_-1&M\ar[r]\ar@{ >->}[d]_-1&0\ar[r]\ar@{ >->}[d]&\cdots\\
\;C:&\cdots\ar@<0.45ex>[r]\ar@<-0.45ex>[r]&0\ar@<0.45ex>[r]\ar@<-0.45ex>[r]&M
\ar@<0.45ex>[r]^-{1+1}\ar@<-0.45ex>[r]_-1&M\ar@<0.45ex>[r]^-1\ar@<-0.45ex>[r]_-1&M
\ar@<0.45ex>[r]\ar@<-0.45ex>[r]&0\ar@<0.45ex>[r]\ar@<-0.45ex>[r]&\cdots}$$
in which $M$ is an abelian monoid and $1=1_M$. Clearly,
$\sg=\linebreak( \dots,0,1,1,1,0,\dots)$ is a morphism of chain
complexes (see 1.14 and 1.11).  One can easi\-ly see that the long
 homology sequence associated to  this diagram coincides  with the sequence
\begin{gather*}
\xymatrix{\cdots\ar[r]&0\ar[r]&E(M)\ar[r]&0\ar[r]&0\ar[r]&0\ar[r]&}\\
\xymatrix{{}\ar[r]&0\ar[r]&M\ar[r]^-k&M'\ar[r]&0\ar[r]&\cdots}
\end{gather*}
where $E(M)$ is the monoid of all idempotents of $M$, $M'$ denotes the largest cancellative
 homomorphic image of $M$, and $k$ is the canonical homomorphism.
($M'=M/\!\!\sim$, \  $m_1\sim m_2$, $m_1,m_2\in M\Leftrightarrow m_1+m=m_2+m$, $m\in M$.
$\cl_\sim(m_1)+\cl_\sim(m_2)=\cl_\sim(m_1+m_2)$, $k(m)=\cl_\sim(m)$.) Let $\dl_1^+=1+1$.
 Hence $H_1(C)=E(M)$, $H_{-1}(B)=M$, $H_{-1}(C)=M'$ and $H_{-1}(\sg)=k$. If $E(M)\neq 0$,
  then this sequence is not exact  at $H_1(C)$ as well as at $H_{-1}(B)$.
\vskip+2mm

{\bf 2.7. Example.} Let $f=\{f_n\}:X=\{X_n,\pa_n^+,\pa_n^-\}\lra X'=\linebreak
\{X'_n,{\pa\,'}_n^+,{\pa_n^{\,'}}^-\}$ be a morphism of chain complexes. The {\em mapping
 cone} of $f$ is the chain complex
\begin{gather*}
C_f=\big\{(C_f)_n,d_n^+,d_n^-\big\},\quad(C_f)_n=X_{n-1}\oplus X'_n,\\
d_n^+(x,x')\!=\!\big(\pa_{n-1}^-(x),{\pa_n^{\,'}}^+(x')\!+\!f_{n-1}(x)\big),\quad
d_n^-(x,x')\!=\!\big(\pa_{n-1}^+(x),{\pa_n^{\,'}}^-(x')\big).
\end{gather*}
There is a Schreier short exact sequence of chain complexes and
their $\pm$-morphisms {\small$$\hskip-2mm\xymatrix{\ \ \ \ \ \ \ \
\  X':\;\cdots\ar@{ ->>}[dd]^-{i_{{}_f}}\ar@<0.45ex>[r]
\ar@<-0.45ex>[r]&X'_{n+1}
\ar@<0.45ex>[r]^-{{\pa^{\,'^+}_{n+1}}}\ar@<-0.45ex>[r]_-{{\pa^{\,'^-}_{n+1}}}
\ar@{ ->>}[dd]^-{(i_{{}_f})_{n+1}}
&X'_n\ar@<0.45ex>[r]^-{{\pa_n^{\,'}}^+}\ar@<-0.45ex>[r]_-{{\pa_n^{\,'}}^-}\ar@{
->>}
[dd]^-{(i_{{}_f})_n}&X'_{n-1}\ar@<0.45ex>[r]\ar@<-0.45ex>[r]\ar@{ ->>}[dd]^-{(i_{{}_f})_{n-1}}&\cdots\\
&&&\\
E_f:\;\;\;C_f:\;\cdots\ar@{
->>}[dd]^-{p_{{}_f}}\ar@<0.45ex>[r]\ar@<-0.45ex>[r]&X_n\!\oplus\!
 X'_{n+1}\ar@<0.45ex>[r]^{d_{n+1}^+}\ar@<-0.45ex>[r]_-{d_{n+1}^-}\ar@{ ->>}[dd]^-{(p_{{}_f})_{n+1}}
&X_{n-1}\!\oplus\!
X'_n\ar@<0.45ex>[r]^-{d_n^+}\ar@<-0.45ex>[r]_-{d_n^-}\ar@{ ->>}
[dd]^-{(p_{{}_f})_n}&X_{n-2}\!\oplus\!
X'_{n-1}\ar@<0.45ex>[r]\ar@<-0.45ex>[r]\ar@{ ->>}[dd]^-
{(p_{{}_f})_{n-1}}&\cdots\\
&&&\\
\;\;\;\;\;\;X[-1]:\cdots\ar@<0.45ex>[r]\ar@<-0.45ex>[r]
&X_n\ar@<0.45ex>[r]^-{\pa_n^-}\ar@<-0.45ex>[r]_-{\pa_n^+}
&X_{n-1}\ar@<0.45ex>[r]^-{\pa_{n-1}^-}\ar@<-0.45ex>[r]_-{\pa_{n-1}^+}
&X_{n-2}\ar@<0.45ex>[r]\ar@<-0.45ex>[r]&\cdots}$$} where
$(i_{{}_f})_n$ sends $x'$ to $(0,x')$, and $(p_{{}_f})_n$ sends
$(x,x')$ to $x$. An element $(x,x')$ of $X_{n-1}\oplus X'_n$ is a
representative of $(E_f)_n$ if and only if
 $x'\in U(X'_n)$, where $U(X'_n)$ denotes the maximal $\Lb$-submodule of $X'_n$. Therefore
 each $d_n^-$ obviously preserves representatives. Assume that $X'$ is a chain complex of
  cancellative $\Lb$-semimodules. Then, by Proposition 2.1, we have the long homology
  sequence
\begin{gather*}
\xymatrix{H(E_f):\ \ \cdots\ \to H_n(X')\ar[r]^-{H_n(i_{{}_f})}
&H_n(C_f)\ar[r]^-{H_n(p_{{}_f})}&H_n(X[-1])\ar[r]^-{\pa_n(E_f)}&}\\
\;\;\;\xymatrix{\ar[r]^-{\pa_n(E_f)}&H_{n-1}(X')\ar[rr]^-{H_{n-1}(i_{{}_f})}
&&H_{n-1}(C_f)\ar[r]&\cdots}
\end{gather*}
associated to $E_f$. One can easily see that in fact $H_n(X[-1])=H_{n-1}(X)$, and
 $\pa_n(E_f)=H_{n-1}(f)$. Furthermore, Theorems 2.4 and 2.5 together with 1.7 imply
  the following
\vskip+2mm

{\bf Corollary.} {\em Let
$f=\{f_n\}:X\!=\!\{X_n,\pa_n^+,\pa_n^-\}\lra X'=\{X'_n,
\pa_n^{'\,^+},\pa_n^{'\,^-}\}$ be a morphism of chain complexes
and suppose that one of the following holds:
\begin{enumerate}
\item[(i)] For each $n$, $X'_n$ is a cancellative
$\Lb$-semimodule,  $X_n$ a $\Lb$-module, and
$H_n(p_{{}_f})(H_n(C_f))=H_n(p_{{}_{K(f)}})(H_n(C_{{}_{K(f)}}))$.

\item[(ii)] For each $n$, $X'_n$ is a $\Lb$-module and $H_n(f)(H_n(X))=\linebreak H_n(K(f))(
H_n(K(X)))$.
\end{enumerate}
Then $H(E_f)$ is exact everywhere, and $H_n(i_{{}_f})$,
$H_n(p_{{}_f})$ and $\pa_n(E_f)(=H_{n-1}(f))$ are
 normal $\Lb$-homomorphisms for all $n$ $($see $2.3)$.}
\vskip+2mm

\begin{proof}
Suppose (i) holds. Since $K$ commutes with mapping cones, it
follows easily that
$H_n(p_{{}_f})(H_n(C_f))=H_n(K(p_{{}_f}))(H_n(K(C_f)))$ for all
$n$. Therefore, thinking of $p_{{}_f}=\{(p_{{}_f})_n\}$ as a
morphism from $C_f$ to $\{X_{n-1},\pa_{n-1}^--\pa_{n-1}^+\}$, we
conclude, by 2.4, that $H(E_f)$ is exact everywhere and
$H_n(i_{{}_f})$ is normal. By 1.7, $H_n(p_{{}_f})$ and
$\pa_n(E_f)$ are
 also normal. When (ii) holds, the assertion is clear since $H_{n-1}(f)=\pa_n(E_f)=
 \pa_n(K(E_f))H_n(k_{X[-1]})$ and $H_{n-1}(K(f))=\pa_n(E_{K(f)})=
 \pa_n(K(E_f))$ (see 2.5 and 1.7).
\end{proof}

In subsequent papers we shall give applications of 2.4, 2.5 and 2.7.

\

\

{\small Alex Patchkoria

A. Razmadze Mathematical Institute

1, M. Aleksidze St., Tbilisi 0193

Georgia

Email: \texttt{lexo@rmi.acnet.ge}}

\end{document}